\documentclass[11pt,a4paper]{amsart}

\usepackage[latin1]{inputenc}
\usepackage{amsmath}
\usepackage{amsfonts}
\usepackage{amssymb}
\usepackage{amsthm}
\usepackage{hyperref}
\usepackage{anysize}
\usepackage{enumitem}
\usepackage[all,cmtip]{xy}

\marginsize{3.0cm}{3.0cm}{3.0cm}{3.0cm}

\newtheorem{theorem}{Theorem}[section]
\newtheorem{lemma}[theorem]{Lemma}
\newtheorem{proposition}[theorem]{Proposition}
\newtheorem{corollary}[theorem]{Corollary}
\newtheorem{conjecture}[theorem]{Conjecture}

\theoremstyle{definition}

\newtheorem{examples}[theorem]{Examples}

\theoremstyle{remark}
\newtheorem{remark}[theorem]{Remark}

\newcommand{\coloneqq}{\mathrel{\mathop:}=}

\newcommand{\puresymbol}[1]{{\lbrace #1 \rbrace}}
\newcommand{\lpuresymbol}[2]{{\lbrace #1, \hdots , #2 \rbrace}}
\newcommand{\logsymbol}[1]{{\frac{d#1}{#1}}} 
\newcommand{\llogsymbol}[2]{{\frac{d#1}{#1} \wedge \hdots \wedge \frac{d#2}{#2}}}

\setenumerate{label={\normalfont(\arabic*)}}

\numberwithin{equation}{section}
\begin{document}

\title[Milnor $K$-groups and function fields of hypersurfaces in characteristic $p$]{Milnor $K$-groups and function fields of hypersurfaces in positive characteristic}
\author{Stephen Scully}
\address{Max-Planck-Institut f\"{u}r Mathematik, Vivatsgasse 7, 53111 Bonn, Germany.}
\email{sscully@mpim-bonn.mpg.de}


\begin{abstract} Let $X$ be an integral affine or projective hypersurface over a field $F$ of characteristic $p>0$, and let $F(X)$ denote its function field. In a recent article, Dolphin and Hoffmann obtained an explicit description of the kernel of the natural restriction homomorphism between the rings of absolute K\"{a}hler differentials of $F$ and $F(X)$, respectively. In this note, we examine the possibility of deriving an analogous result for mod-$p$ Milnor $K$-theory using the Bloch-Gabber-Kato theorem. \end{abstract}
\maketitle

\section{Introduction} 

Let $F$ be a field and let \textsf{\textsl{Fields}}$_F$ denote the category of fields containing $F$ as a subfield (with morphisms given by $F$-linear field homomorphisms). For every non-negative integer $n$ and prime integer $p$, let  $k_{n,p} \colon \textsf{\textsl{Fields}}_F \rightarrow \textsf{\textsl{AbelianGroups}}$ be the functor which assigns to a given overfield $K$ of $F$ the group $k_{n,p}(K) \coloneqq K_n^M(K)\otimes \mathbb{Z}/p$, where $K_n^M(K)$ denotes the $n^{th}$ Milnor $K$-group of $K$ (see \cite{MilnorInventiones} or \S \ref{MilnorKtheory} below). Given an inclusion $K \subseteq L$ of fields containing $F$, we write $k_{n,p}(L/K)$ for the kernel of the associated group homomorphism $k_{n,p}(K) \rightarrow k_{n,p}(L)$. If $X$ is an (integral) $F$-variety with function field $F(X)$, then the groups $k_{n,p}(F(X)/F)$ comprise an important collection of birational invariants of $X$. On the one hand, these invariants encode non-trivial information concerning the structure of $X$ as a geometric object, constituting particular components of the localisation sequences
\begin{equation*} k_{n,p}(F) \rightarrow k_{n,p}\big(F(X)\big) \rightarrow \bigoplus_{x \in X^{(1)}} k_{n-1.p}\big(F(x)\big) \rightarrow \bigoplus_{x \in X^{(2)}} k_{n-2,p}\big(F(x)\big) \rightarrow \hdots \end{equation*}
(see \cite{KatoMilnorKTheoryZeroCycles}); on the other, their study (for suitable classes of varieties $X$) unifies several fundamental problems relating to the classification of certain classical algebraic objects defined over fields, notably in the context of the following well-known examples:
\begin{examples} \label{Examples1} Let $K$ be a field.
\begin{enumerate}[leftmargin=*] \item If $\mathrm{char}(K) \neq p$, then there exists, for any $n \geq 0$, a natural isomorphism between $k_{n,p}(K)$ and the Galois cohomology group $H^n(K, \mu_p^{\otimes n})$ by a celebrated recent result of Rost and Voevodsky (formerly the (motivic) ``Bloch-Kato conjecture'', cf. \cite{VoevodskyMotivicCohomology}).
\item If $\mathrm{char}(K) \neq p$ and $K$ contains a primitive $p^{th}$ root of unity, then, as a special case of (1), there exists a natural isomorphism $k_{2,p}(K) \simeq {_p\mathrm{Br}(K)}$, where $_p\mathrm{Br}(K)$ denotes the $p$-torsion part of the Brauer group of $K$ (this result was originally proved by Merkurjev and Suslin in \cite{MerkurjevSuslinNormResidue}).
\item For any $n \geq 0$, there exists a natural isomorphism between $k_{n,2}(K)$ and the degree-$n$ part of the graded Witt ring of symmetric bilinear forms over $K$ by celebrated results of Kato (in the case where $\mathrm{char}(K) = 2$, \cite{KatoInventiones}) and Orlov-Vishik-Voevodsky (in the case where $\mathrm{char}(K) \neq 2$, \cite{OrlovVishikVoevodskyAnnals}). \end{enumerate} \end{examples}

In the case where $\mathrm{char}(F) \neq p$, determining the groups $k_{n,p}(F(X)/F)$ for a non-unirational $F$-variety $X$ is an extremely difficult task in general (for unirational $X$ we have $k_{n,p}(F(X)/F) =0$ for all $n\geq 0$ by \cite[Theorem 2.3]{MilnorInventiones}). For instance, in the context of Example \ref{Examples1} (3), the problem of computing the groups $k_{n,2}(F(X)/F)$ in the case where $\mathrm{char}(F) \neq 2$ and $X$ is an (anisotropic) $F$-quadric already amounts to one of the major outstanding problems in the theory of quadratic forms over general fields. In this case, we have the following explicit conjecture due to Vishik (cf. \cite[Question 1]{VishikKernelsMilnorsKTheory}):

\begin{conjecture} \label{Vishiksconjecture} Let $X$ be an affine or projective quadric over a field $F$ of characteristic $\neq 2$. Then, for any $n \geq 0$, $k_{n,2}(F(X)/F)$ is generated by (the classes of) symbols. \end{conjecture}

By a deep result due to Orlov, Vishik and Voevodsky, this conjecture is known to be valid in the case where $n < \mathrm{log}_2(\mathrm{dim}\;X +2)$ (in fact, Orlov et al. show that $k_{n,2}(F(X)/F) = 0$ in this case, cf. \cite[Theorem 4.2]{OrlovVishikVoevodskyAnnals}). A positive solution to Conjecture \ref{Vishiksconjecture} in the next open case (where $n = \mathrm{log}_2[2\mathrm{dim}\;X + 2]$) would have significant ramifications for the theory of quadratic forms over fields of characteristic different from 2; most notably, a conjectural characterisation of the important class of so-called \emph{Pfister neighbours} in terms of ``motivic'' data would follow immediately (see \cite{VishikKernelsMilnorsKTheory},\cite{IzhboldinVishikMaximalSplitting} for further discussion).

By contrast, the picture simplifies considerably in the case where $\mathrm{char}(F) = p$. Indeed, suppose that we are in this setting, and let $\Omega^n \colon \textsf{\textsl{Fields}}_F \rightarrow \textsf{\textsl{AbelianGroups}}$ be the functor which assigns to a given overfield $K$ of $F$ the group $\Omega^n_K$ of absolute differential $n$-forms of $K$ (see \cite[\S A.8]{GilleSzamuelyGaloisCohomology} or \S \ref{Kahlerdifferentials} below). Then, by a celebrated theorem of Bloch, Gabber and Kato (cf. \cite[Theorem 2.1]{BlochKatopadicetalecohomology} or Theorem \ref{BlochGabberKatoTheorem} below), the assignment $\lpuresymbol{a_1}{a_n} \mapsto \llogsymbol{a_1}{a_n}$ gives rise, for any $n \geq 0$, to an injective morphism of functors $h \colon k_{n,p} \rightarrow \Omega^n$. In particular, for any $F$-variety $X$ and any $n \geq 0$, we have an identification
\begin{equation} \label{eq1.1} k_{n,p}(F(X)/F) \simeq h_F\big(k_{n,p}(F)\big)\!\cap \Omega^n(F(X)/F), \end{equation}
where $\Omega^n(F(X)/F)$ denotes the kernel of the group homomorphism $\Omega^n_F \rightarrow \Omega^n_{F(X)}$ induced by the (structural) inclusion of $F$ into $F(X)$. The significance of these identifications lies in the fact that the groups $\Omega^n(F(X)/F)$ (being, in reality, linear spaces) are more readily accessible than those which we ultimately aim to compute. For example, one of the basic results of the elementary theory of K\"{a}hler differentials is the uniform triviality of the groups $\Omega^n(F(X)/F)$ in the case where $X$ is \emph{generically smooth} (i.e. where $F(X)$ is separable over $F$, see Proposition \ref{PROPgenericsmoothness} or Lemma \ref{LEMSeparabletriviality} below). Thus, in view of \ref{eq1.1}, we immediately see that $k_{n,p}(F(X)/F) = 0$ for all $n \geq 0$ in this case. Corroborating this point further, Dolphin and Hoffmann (see \cite[Theorem 8.5]{DolphinHoffmannDifferentialandBilinearForms} or \S \ref{Kahlerdifferentialsandfunctionfieldsofhypersurfaces} below) recently obtained an explicit description of the groups $\Omega^n(F(X)/F)$ for an \emph{arbitrary} integral $F$-\emph{hypersurface} $X$. In the \emph{nowhere-smooth} case (the only case of interest in light of the above remarks), this description is made in terms of another birational invariant of $X$ known as its \emph{norm field}, which may be interpreted as the smallest algebraic overfield $N(X)$ of $F$ such that the regular locus of the scheme $X \otimes_F N(X)$ is empty (we should remark here that the invariant defined in \cite{DolphinHoffmannDifferentialandBilinearForms} is actually the $p^{th}$-power subfield of our $N(X)$; since $p = \mathrm{char}(F)$, however, any distinction to be drawn between the two is immaterial, see \S \ref{NormField} below). More explicitly, if $X$ is defined as a hypersurface in the ambient affine (resp. projective) space by the vanishing of a given polynomial (resp. homogeneous polynomial) $f$, then $N(X)$ is constructed by adjoining to $F$ the $p^{th}$ roots of all possible ratios of the (non-zero) coefficients of $f$ (note that $F$ is necessarily imperfect); in particular, $N(X)$ is a height-one finite purely inseparable extension of $F$. With this terminology in place, the main result of \cite{HoffmannDiagonalFormsofDegreepinCharp} asserts that, for any nowhere-smooth integral $F$-hypersurface $X$ and any $n \geq 0$, $\Omega^n(F(X)/F)$ consists of those elements of $\Omega_F^n$ which are divisible by the $m$-fold exterior product $\llogsymbol{b_1}{b_m}$, where $b_1,\hdots,b_m$ is any minimal generating set for the extension $N(X)^p/F^p$ (see also Theorem \ref{DolphinHoffmanntheorem} and Remark \ref{Explicitremark} below). This computation was successfully used (in conjunction with the celebrated results of Kato (\cite{KatoInventiones}) on 2-primary analogues of the Milnor conjectures) in \emph{loc. cit.} to study the splitting behaviour of symmetric bilinear forms under scalar extension to function fields of hypersurfaces in characteristic 2. In the context of Milnor $K$-theory, it leads naturally to the following conjecture (communicated to the author by Hoffmann):

\begin{conjecture} \label{Mainconjecture} Let $X$ be a nowhere-smooth integral affine or projective hypersurface over an imperfect field $F$ of characteristic $p>0$. Then, for any $n \geq 0$, $k_{n,p}(F(X)/F)$ is generated by (the classes of) those symbols $\lpuresymbol{a_1}{a_n}$ for which $N(X) \subseteq F(\sqrt[p]{a_1},\hdots,\sqrt[p]{a_n})$. \end{conjecture}

The difficulty in establishing assertions of this type via the identifications of \ref{eq1.1} is that, given a non-trivial $F$-linear subspace $V \subseteq \Omega^n_F$, there is no general method by which one can derive an explicit presentation of the pre-image of $V$ under the Bloch-Gabber-Kato embedding $h_F \colon k_{n,p}(F) \rightarrow \Omega_F^n$ in terms of symbols. In the case at hand, for example, the missing link between Conjecture \ref{Mainconjecture} and the main result of \cite{DolphinHoffmannDifferentialandBilinearForms} is given by: 

\begin{conjecture} \label{Secondaryconjecture} Let $F$ be a field of characteristic $p>0$, let $b_1, \hdots, b_m \in F$ be such that $[F^p(b_1,\hdots,b_m):F^p] = p^m$ and let $V = \Omega_F^{n-m} \wedge \big(\llogsymbol{b_1}{b_m}\big) \subseteq \Omega_F^n$ for some $n \geq m$. Then $h_F^{-1}(V) \subseteq k_{n,p}(F)$ is generated by (the classes of) those symbols $\lbrace a_1,\hdots,a_n \rbrace$ for which $F(\sqrt[p]{b_1},\hdots, \sqrt[p]{b_m}) \subseteq F(\sqrt[p]{a_1},\hdots,\sqrt[p]{a_n})$. \end{conjecture}

In the present article, we will show that Conjecture \ref{Secondaryconjecture} is valid in the special case where $n = m$ (see Proposition \ref{PROPproofofsecondaryconjecture} below). Together with our previous remarks, this yields:

\begin{theorem} \label{Maintheorem} Let $X$ be an integral affine or projective hypersurface over a field $F$ of characteristic $p>0$.
\begin{enumerate} \item If $X$ is generically smooth, then $k_{n,p}(F(X)/F) = 0$ for all $n \geq 0$.
\item If $X$ is nowhere smooth, then $k_{n,p}(F(X)/F) = 0$ for all $0 \leq n < \mathrm{log}_p[N(X):F]$.
\item If $X$ is nowhere smooth and $m = \mathrm{log}_p[N(X):F]$, then $k_{m,p}(F(X)/F)$ is generated by (the classes of) those symbols $\lpuresymbol{a_1}{a_m}$ for which $N(X) = F(\sqrt[p]{a_1},\hdots,\sqrt[p]{a_m})$.
\end{enumerate} \end{theorem}

In particular, we see that the kind of behaviour which is expected to prevail for quadric hypersurfaces over fields of characteristic different from 2 (explicitly stated in Conjecture \ref{Vishiksconjecture} above and the subsequent remarks) also prevails in this setting.\\

\noindent {\bf Notation and terminology.} If $F$ is a field, then $\overline{F}$ will denote a fixed algebraic closure of $F$, while \textsf{\textsl{Fields}}$_F$ and \textsf{\textsl{VectorSpaces}}$_F$ will denote the category of fields containing $F$ as a subfield (with morphisms given by $F$-linear field homomorphisms) and the category of $F$-vector spaces, respectively. The word \emph{scheme} means an scheme of finite type over a field. By a \emph{variety}, we then mean an integral scheme. If $X$ (resp. $f$) is a scheme over a field $F$ (resp. a polynomial in $n$-algebraically independent variables $T_1,\hdots,T_n$ over $F$) and $L$ is any overfield of $F$, then $X \otimes_F L$ (resp. $f \otimes_F L$) will denote the scheme $X \times_{\mathrm{Spec}\;F} \mathrm{Spec}\;L$ (resp. the image of $f$ under the canonical inclusion $F[T_1,\hdots,T_n] \subseteq L[T_1,\hdots,T_n]$).

\section{Preliminaries} \label{Preliminaries}

The purpose of this section is to recall some basic facts concerning Milnor $K$-groups and modules of (higher) differential forms. Having no need to work in greater generality, we limit our considerations here to the restriction of these functors to categories of \emph{fields}.

\subsection{Milnor $K$-theory} \label{MilnorKtheory} Let $K$ be a field. The \emph{Milnor ring} of $K$ (first introduced in \cite{MilnorInventiones}) is the $\mathbb{Z}$-graded ring $K_*^M(K)$ defined as the quotient of the tensor ring $T^*_{\mathbb{Z}}(K^*)$ (equipped with its standard $\mathbb{Z}$-grading) by the homogeneous ideal generated by elements of the form $a \otimes (1-a)$ for some $a \in K^* \setminus \lbrace 1 \rbrace$. Its $n^{th}$ graded component, $K_n^M(K)$, is called the $n^{th}$ \emph{Milnor $K$-group} of $K$. For any $a_1,\hdots,a_n \in K^*$, the image of the $n$-fold product $a_1 \otimes \hdots \otimes a_n$ under the canonical projection $T^n_{\mathbb{Z}}(K^*) \rightarrow K_n^M(K)$ is denoted by $\lpuresymbol{a_1}{a_n}$. Elements of the latter type are known as (pure) \emph{symbols}; collectively (over all $n\geq 0$), they form a system of additive generators of the ring $K_*^M(K)$. Given a prime integer $p$, we will write $k_{n,p}(K)$ for the mod-$p$ quotient group $K_n^M(K) \otimes \mathbb{Z}/p$.

Let $L$ be another field, and let $i \colon K \rightarrow L$ be a field homomorphism. Then the induced map $T_{\mathbb{Z}}^*(K^*) \rightarrow T_{\mathbb{Z}}^*(L^*)$ (defined as the unique ring homomorphism extending $i|_{K^*}$) descends to a $\mathbb{Z}$-graded ring homomorphism $K_*^M(K) \rightarrow K_*^M(L)$, thus equipping $K_*^M(L)$ with the structure of a $\mathbb{Z}$-graded $K_*^M(K)$-module. In particular, for any $n \geq 0$ and any prime integer $p$, $i$ gives rise to a \emph{restriction homomorphism} 
\begin{equation*} r_{i} \colon k_{n,p}(K) \rightarrow k_{n,p}(L). \end{equation*}
These restriction maps are compatible with the composition of field homomorphisms, and, in this way, we obtain functors $k_{n,p} \colon \textsf{\textsl{Fields}}_F \rightarrow \textsf{\textsl{AbelianGroups}}$ for any fixed base field $F$. If $i$ is the inclusion of $K$ as a subfield of $L$, then we simply write $r_{K \subseteq L}$ instead of $r_i$; in this case, we also write $k_{n,p}(L/K)$ for the kernel of $r_i$ (in degree $n$).

In the special case where $i$ is \emph{finite} (i.e. where $i$ equips $L$ with the structure of a \emph{finite-dimensional} $K$-vector space), there also exist \emph{corestriction homomorphisms}
\begin{equation*} c_{i} \colon k_{n,p}(L) \rightarrow k_{n,p}(K), \end{equation*}
which are again compatible with the composition of field homomorphisms. In the case where $n=1$, the map $c_i$ is simply given by the standard (mod-$p$) \emph{norm homomorphism}; for larger $n$, however, the construction is rather more delicate (see \cite{BassTateMilnorRingGlobalField},\cite{KatoLocalClassFieldTheoryII} or \cite[Ch. 7]{GilleSzamuelyGaloisCohomology} for full details). Together, the restriction and corestriction homomorphisms satisfy the \emph{projection formula}: for any $n \geq 0$ and prime integer $p$, the composite map $c_i \circ r_i \in \mathrm{End}\big(k_{n,p}(K)\big)$ amounts to multiplication by $\mathrm{deg}(i) = [L:i(K)]$. As before, in the case where $i$ is the inclusion of $K$ as a subfield of $L$, we simply write $c_{K \subseteq L}$ instead of $c_i$.

\subsection{K\"{a}hler differentials} \label{Kahlerdifferentials} Let $k \subseteq K$ be an inclusion of fields. The \emph{space of differential 1-forms of $K$ relative to $k$} is defined as the (essentially) unique object $\Omega_{K|k}^1$ which represents the functor $\mathrm{Der}_k(K,-) \colon \textsf{\textsl{VectorSpaces}}_K \rightarrow \textsf{\textsl{Sets}}$ of $k$-linear derivations of $K$ (see \cite[Ch. 9]{MatsumuraCommutativeRingTheory}) in the category of $K$-vector spaces. Explicitly, it may be presented as the quotient $V/W$, where $V$ is the free $K$-vector space on the set $\lbrace (a)\;|\;a \in K \rbrace$ and $W$ is the $K$-linear subspace of $V$ generated by elements of the form $(\mu)$, $(a + b) - (a) - (b)$ or $(ab) - b(a) - a(b)$ for some $\mu \in k$ or $a,b \in K$. The universal $k$-linear derivation $d \colon K \rightarrow \Omega_{K|k}^1$ is then defined by the assignment $a \mapsto da \coloneqq (a) \pmod{W}$. By setting
\begin{equation} \label{differentialdefinition} d(ada_1 \wedge \hdots \wedge da_n) = da \wedge da_1 \wedge \hdots \wedge da_n \end{equation}
for all $n \geq 0$ and all $a, a_1, \hdots, a_n \in K$, the latter map may be extended to a unique degree-one graded $k$-algebra endomorphism $d$ of the exterior algebra $\Omega_{K|k}^* \coloneqq \bigwedge_K^*(\Omega_{K|k}^1)$ (equipped with its standard $\mathbb{Z}$-grading) such that $d^2=0$ and such that the formula
\begin{equation*} d(\omega_1 \wedge \omega_2) = d(\omega_1) \wedge \omega_2  + (-1)^i \omega_1 \wedge d(\omega_2) \end{equation*}
holds for all $i \geq 0$, $\omega_1 \in \Omega_{K|k}^i$ and $\omega_2 \in \Omega_{K|k}^*$. In particular, letting $\Omega_{K|k}^n$ denote the degree-$n$ part of the ring $\Omega_{K|k}^*$, we obtain in this way a complex
\begin{equation*} (\Omega_{K|k}^\bullet,d) = K \xrightarrow{d} \Omega_{K|k}^1 \xrightarrow{d} \Omega_{K|k}^2 \xrightarrow{d} \hdots \end{equation*}
of $k$-vector spaces known as the \emph{de Rham complex of $K$ relative to $k$}. For each $n \geq 1$, the $K$-vector space $\Omega_{K|k}^n$ is called the \emph{space of differential $n$-forms of $K$ relative to $k$}. Its $k$-linear subspace consisting of the $n$-coboundaries of $(\Omega_{K|k}^\bullet,d)$ (i.e. the elements in image of the differential $d \colon \Omega_{K|k}^{n-1} \rightarrow \Omega_{K|k}^n$) will be denoted by $B_{K|k}^n$. For any $n \geq 0$, the $n^{th}$ cohomology set of the complex $(\Omega_{K|k}^\bullet,d)$ (which is naturally equipped with the structure of a $k$-vector space) will be denoted by $H^n(\Omega^\bullet_{K|k})$. If $\omega \in \Omega_{K|k}^n$ lies in the kernel of the differential $d$, then we write $[\omega]$ for its image in $H^n(\Omega^\bullet_{K|k})$. In the special case where $k$ is the prime subfield of $K$, we simply write $\Omega_K^n$ \big(resp. $B_K^n$, $H^n(\Omega_K^\bullet)$\big) instead of $\Omega_{K|k}^n$ \big(resp. $B_{K|k}^n$, $H^n(\Omega_{K|k}^\bullet)$\big), and call it the \emph{space of absolute differential $n$-forms of $K$}.

Let $\ell \subseteq L$ be another inclusion of fields, and let $i \colon K \rightarrow L$ and $j \colon k \rightarrow \ell$ be a pair of field homomorphisms such that $i|_k$ is equal to the composition of $j$ and the inclusion of $\ell$ into $L$ (note that although the action of $j$ is completely determined by $i$, we want to emphasise the target $\ell$ of $j$ in what follows). Consider the $k$-linear derivation $D_{i,j} \colon K \rightarrow \Omega_{L|\ell}^1$ defined by setting $D_{i,j}(a) = di(a)$ for all $a \in K$. By the universal property of $\Omega_{K|k}^1$, $D_{i,j}$ corresponds to a unique $K$-linear homomorphism $\Omega_{K|k}^1 \rightarrow \Omega_{L|\ell}^1$. Taking $n^{th}$ exterior powers, we therefore obtain a $K$-linear \emph{restriction homomorphism}
\begin{equation*} r_{i,j} \colon \Omega_{K|k}^n \rightarrow \Omega_{L|\ell}^n \end{equation*}
for all $n \geq 0$. More explicitly, $r_{i,j}$ is defined by setting $r_{i,j}(ada_1 \wedge \hdots \wedge da_n) = i(a) di(a_1) \wedge \hdots \wedge di(a_n)$ for all $a,a_1,\hdots,a_n \in K$ and then extending linearly to all of $\Omega_{K|k}^n$. If $j$ is the inclusion of $k$ as a subfield of $\ell$ (resp. $i$ is the inclusion of $K$ as a subfield of $L$), we simply write $r_{i|k \subseteq \ell}$ (resp. $r_{K \subseteq L|j}$) instead of $r_{i,j}$. In the special case where $\ell = k$ and $j$ is the identity (resp. $L=K$ and $i$ is the identity), we simplify the notation further, denoting the appropriate restriction homomorphism by $r_i$ (resp. $r_{|j}$). Finally, if $i$ is the inclusion of $K$ as a subfield of $L$ and $\ell = k = k_0$, where $k_0$ denotes the prime subfield of $K$, then we write $\Omega^n(L/K)$ for the kernel of the map $r_{K \subseteq L}$ (in degree $n$). The maps $r_{i,j}$ are readily seen to respect the composition of compatible pairs of field homomorphisms. In particular, for any fixed base field $F$ and any $n \geq 0$, we obtain a functor $\Omega^n \colon \textsf{\textsl{Fields}}_F \rightarrow \textsf{\textsl{AbelianGroups}}$ which assigns to a given overfield $K$ of $F$ the underlying group of the $K$-vector space $\Omega_K^n$. The restriction homomorphisms also commute with the differentials of the de Rham complex, as one may readily verify from the definition and \eqref{differentialdefinition}.

A totally ordered collection $\lbrace a_i \rbrace_{i \in I}$ of elements of $K$ is said to be a \emph{differential basis of $K$ over $k$} if, for every $n \geq 0$, the set
\begin{equation*}  \mathfrak{B}_n = \lbrace \llogsymbol{a_{s(1)}}{a_{s(n)}}\;|\;s \colon \lbrace 1,\hdots,n \rbrace \rightarrow I \text{ is a strictly increasing function} \rbrace \end{equation*}
constitutes a basis of the $K$-vector space $\Omega_{K|k}^n$ (in the case where $n=0$, we have $\Omega_{K|k}^0 = K$, and the set $\mathfrak{B}_0$ should be interpreted as the one-element set $\lbrace 1 \rbrace$). It is not difficult to see (cf. \cite[Theorem 26.5]{MatsumuraCommutativeRingTheory}) that this notion coincides with that of a \emph{transcendence basis} of $K$ over $k$ in the case where $\mathrm{char}(k) = 0$, and that of a \emph{$p$-basis} of $K$ over $k$ in the case where $\mathrm{char}(k) = p>0$ (recall that $\lbrace a_i \rbrace_{i \in I}$ is a \emph{$p$-basis of $K$ over $k$} if $K = K^p(k)(a_i\;|\; i \in I)$ and $[K^p(k)(a_{s(1)}, \hdots, a_{s(n)}):K^p(k)] = p^n$ for every strictly increasing function $s \colon \lbrace 1,\hdots,n \rbrace \rightarrow I$). The following proposition therefore follows immediately from standard results on the extension of transcendence bases and $p$-bases (in the positive-characteristic case, we refer in particular to the well-known separability criterion of S. MacLane, see \cite{MacLaneModularFields}):

\begin{proposition}[{cf. \cite[Theorem 26.6]{MatsumuraCommutativeRingTheory}}] \label{PROPseparability} Let $K \subseteq L$ be an inclusion of fields. Then $L$ is separable over $K$ if and only if, for any subfield $k \subseteq K$ and any $n \geq 0$, the $L$-linear homomorphism $\alpha \colon \Omega_{K|k}^n \otimes_K L \rightarrow \Omega_{L|k}^n$ induced by $r_{K\subseteq L}$ is injective. \end{proposition}

\begin{remark} Recall here that $L$ is \emph{separable} over $K$ if the ring $L \otimes_K \overline{K}$ is reduced. \end{remark}

In the special case where $L$ is algebraic over $K$, one can say more:

\begin{proposition} \label{PROPseparablealgebraic} Let $k \subseteq K \subseteq L$ be a tower of fields, and let $\ell$ be a subfield of $L$ containing $k$. Suppose that $L$ is separable and algebraic over $K$, and that $k$ is separable and algebraic over $\ell$. Then, for any $n \geq 0$, the restriction map $r_{K \subseteq L|k \subseteq \ell} \colon \Omega_{K|k}^n \rightarrow \Omega_{L|\ell}^n$ gives rise to an isomorphism of $L$-vector spaces $\beta \colon \Omega_{K|k}^n \otimes_K L \xrightarrow{\sim} \Omega_{L|\ell}^n$. \end{proposition}

For lack of a precise reference, we provide the details:

\begin{proof} Let $\alpha \colon \Omega_{K|k}^n \otimes_K L \rightarrow \Omega_{L|k}^n$ be the $L$-linear injection of Proposition \ref{PROPseparability}. The proposed isomorphism $\beta$ evidently decomposes as
\begin{equation*} \Omega_{K|k}^n \otimes_K L \xrightarrow{\alpha} \Omega_{L|k}^n \xrightarrow{r_{|k \subseteq \ell}} \Omega_{L|\ell}^n. \end{equation*}
We claim that both $\alpha$ and $r_{|k \subset \ell}$ are $L$-linear isomorphisms. It is sufficient to treat the case where $n=1$ (the general case follows by taking exterior powers). Then, according to the relative cotangent sequence for differential 1-forms (\cite[Theorem 25.1]{MatsumuraCommutativeRingTheory}), the cokernel of $\alpha$ is given by $\Omega_{L|K}^1$. But if $a$ is any element of $L$ with minimal polynomial $m_a$ over $K$, then the relations $0 = d\big(m_a(a)\big) = m_a'(a)da$ (here $m_a'$ denotes the formal derivative of $m_a$) imply (in view of the separability of $a$) that $da = 0$ in $\Omega_{L|K}^1$. It follows that $\Omega_{L|K}^1 = 0$, and so $\alpha$ is injective. In exactly the same way, we see that $\Omega_{\ell|k}^1 = 0$, and hence $r_{|k \subset \ell}$ is injective. As the latter map is clearly surjective, the result follows. \end{proof}

The proof of Proposition \ref{PROPseparablealgebraic} shows that $\Omega_{L|K}^1 =0$ whenever $L$ is separable and algebraic over $K$. It is not hard to see that the converse is also true. More generally, we have:

\begin{proposition}[{see \cite[Theorem 26.10]{MatsumuraCommutativeRingTheory}}] \label{PROPseparabilityfinitelygenerated} Let $K \subseteq L$ be an inclusion of fields with $L$ finitely generated over $K$. Then $L$ is separable over $K$ if and only if $\mathrm{dim}_{L}\Omega_{L|K}^1 = \mathrm{trdeg}_K(L)$. \end{proposition}

Now, let $k \subseteq K \subseteq L$ be a tower of fields, and let $\ell$ be a subfield of $L$ containing $k$. If $L$ is \emph{separable} and \emph{finite} over $K$, and $\ell$ is \emph{separable} and \emph{algebraic} over $k$, then Proposition \ref{PROPseparablealgebraic} shows that the map $\beta \colon \Omega_{K|k}^n \otimes_K L \rightarrow \Omega_{L|\ell}^n$ induced by $r_{K \subseteq L|k \subseteq \ell}$ is an $L$-linear isomorphism for any $n \geq 0$. This enables us to define, for any $n \geq 0$, a \emph{corestriction homomorphism}
\begin{equation*} c_{K \subseteq L|k \subseteq \ell} \colon \Omega_{L|\ell}^n \rightarrow \Omega_{K|k}^n \end{equation*}
as the $K$-linear composition
\begin{equation*} \Omega_{L|\ell}^n \xrightarrow{\beta^{-1}} \Omega_{K|k}^n \otimes_K L \xrightarrow{\mathrm{id} \otimes \mathrm{Tr_{K \subseteq L}}} \Omega_{K|k}^n, \end{equation*}
where $\mathrm{Tr}_{K \subseteq L} \colon L \rightarrow K$ denotes the standard \emph{trace homomorphism} for finite extensions. These maps are easily seen to respect the composition of compatible pairs of field inclusions. Together the restriction and corestriction homomorphisms for finite separable extensions satisfy the \emph{projection formula}: for any $n \geq 0$, the composite map $c_{K \subseteq L|k \subseteq \ell} \circ r_{K \subseteq L|k \subseteq \ell} \in \mathrm{End}_K(\Omega_{K|k}^n)$ amounts to multiplication by $[L:K]$. In the special case where $\ell= k$ (resp. $L=K$), we will again simply write $c_{K\subseteq L}$ (resp. $c_{|k \subseteq \ell}$) instead of $c_{K \subseteq L|k \subseteq \ell}$.

\begin{remark} It is also possible to define corestriction homomorphisms on spaces of higher K\"{a}hler differentials for finite \emph{inseparable} extensions, but the construction is rather more delicate. As the separable transfer is sufficient for our later needs, we refrain from discussing the general case here (see \cite[p. 126]{BlochKatopadicetalecohomology} for further details). \end{remark}

\subsection{A decomposition of the de Rham complex in positive characteristic} \label{Kahlerdifferentialsinpositivecharacteristic}

Let $k \subseteq K$ be an inclusion of fields of characteristic $p>0$ such that $k$ contains $K^p$. Assume that $K$ is finite over $k$, say of degree $p^r$, and let $\lbrace a_i\rbrace_{i \in \lbrace 1,\hdots,r \rbrace}$ be a $p$-basis of $K$ over $k$. For all $0 \leq n \leq r$, let $S_n$ denote the set of all strictly increasing functions $\lbrace 1,\hdots,n \rbrace \rightarrow \lbrace 1, \hdots, r \rbrace$ ($S_0$ should be interpreted here as the one-element set $\lbrace 0 \rbrace$). For any $s \in S_n$, let
\begin{equation*} \omega_s = \llogsymbol{a_{s(1)}}{a_{s(n)}} \in \Omega_{K|k}^n \end{equation*}
(with $\omega_0 = 1 \in K=\Omega_{K|k}^0$ in the case where $n=0$). By the discussion of \S \ref{Kahlerdifferentials}, the set $\mathfrak{B}_n = \lbrace \omega_s\;|\;s\in S_n \rbrace$ is, for any $n \geq 0$, a basis of the $K$-vector space $\Omega_{K|k}^n$. Now, for each multi-index $\alpha = (\alpha_1,\hdots,\alpha_r) \in \lbrace 0,\hdots,p-1 \rbrace^{\times r}$, let $a^\alpha$ denote the product $a_1^{\alpha_1}\cdots a_r^{\alpha_r} \in K$. Since $\lbrace a_i \rbrace_{i \in \lbrace 1,\hdots,r \rbrace}$ is a $p$-basis of $K$ over $k$, the set $\lbrace a^{\alpha}\;|\;\alpha \in \lbrace 1,\hdots,p-1 \rbrace^{\times r} \rbrace$ is a basis of $K$ as a $k$-vector space. For each $n \geq 0$, the set $\lbrace a^{\alpha}\omega_s\;|\;s \in S_n,\alpha \in \lbrace 1,\hdots,p-1 \rbrace^{\times r} \rbrace$ therefore constitutes a $k$-basis of the $K$-vector space $\Omega_{K|k}^n$. Let $\alpha \in \lbrace 0,\hdots,p-1 \rbrace^{\times r}$. For each $0 \leq n \leq r$, let $\Omega_{K|k}^n(\alpha)$ denote the $k$-linear subspace of $\Omega_{K|k}^n$ generated by the set $\lbrace a^{\alpha}\omega_s\;|\;s \in S_n \rbrace$. If $n < r$, then it is easy to see that the $k$-linear differential $d \colon \Omega_{K|k}^n \rightarrow \Omega_{K|k}^{n+1}$ maps $\Omega_{K|k}^n(\alpha)$ to $\Omega_{K|k}^{n+1}(\alpha)$. In this way, we obtain a $k$-linear subcomplex $(\Omega_{K|k}^\bullet(\alpha),d)$ of the de Rham complex $(\Omega_{K|k}^\bullet,d)$. Since the set $\lbrace a^{\alpha}\omega_s\;|\;s \in S_n,\alpha \in \lbrace 1,\hdots,p-1 \rbrace^{\times r} \rbrace$ is a $k$-basis of the space $\Omega_{K|k}^n$ for every $0 \leq n \leq r$, we thus obtain a decomposition
\begin{equation} \label{additivedecompositiondeRham} (\Omega_{K|k}^\bullet,d) = \bigoplus_{\alpha} (\Omega_{K|k}^\bullet(\alpha),d) \end{equation}
of complexes of $k$-vector spaces. In particular, for any $0 \leq n \leq r$, we have a $k$-linear decomposition $H^n(\Omega_{K|k}^\bullet) = \bigoplus_{\alpha}H^n\big(\Omega_{K|k}^\bullet(\alpha)\big)$, where $H^n\big(\Omega_{K|k}^\bullet(\alpha)\big)$ denotes the degree-$n$ cohomology of the complex $\Omega_{K|k}^\bullet(\alpha)$. As the following proposition shows, these decompositions are very useful from a computational perspective (note that $\underline{0}$ denotes here the unique element of $\lbrace 0,\hdots,p-1 \rbrace^{\times r}$ whose entries are all zero):

\begin{proposition}[{cf. \cite[Lemma 1]{KatoGaloisCohomologyCompleteDiscrete}}] \label{PROPCohomologydeRhamcharp} In the above situation:
\begin{enumerate} \item The differentials of the complex $\Omega_{K|k}^\bullet(\underline{0})$ are trivial.
\item If $\alpha \neq \underline{0}$, then $H^n\big(\Omega_{K|k}^\bullet(\alpha)\big) = 0$ for all $0 \leq n \leq r$. \end{enumerate}
In particular, for any $0 \leq n \leq r$, we have
\begin{equation*} H^n(\Omega_{K|k}^\bullet) = \bigoplus_s k \cdot [\omega_s], \end{equation*}
the sum being taken over the set of all $s \in S_n$.
\begin{proof} See \cite[Proposition 9.4.6]{GilleSzamuelyGaloisCohomology}. \end{proof} \end{proposition}

\subsection{The Artin-Schreier homomorphism} \label{ArtinSchreierhomomorphism} Let $k \subseteq K$ be an inclusion of fields of characteristic $p>0$. For any $n \geq 0$, there exists a unique homomorphism
\begin{equation} \label{ArtinSchreier} \wp_{K|k} \colon \Omega_{K|k}^n \rightarrow \Omega_{K|k}^n/B_{K|k}^n \end{equation}
of abelian groups (known as the \emph{Artin-Schreier homomorphism}) satisfying
\begin{equation} \label{ActionofArtinSchreier} \wp_{K|k}(a \llogsymbol{a_1}{a_n}) = (a^p - a)\llogsymbol{a_1}{a_n} \pmod{B_{K|k}^n} \end{equation}
for all $a,a_1,\hdots,a_n \in K$ (the well-definedness of the assignment \ref{ActionofArtinSchreier} is due to Cartier (\cite{CartierQuestionsRationalite}); for a detailed exposition, we refer to \cite[Ch. 9]{GilleSzamuelyGaloisCohomology}). Its kernel will be denoted by $\nu(n)_{K|k}$. If $\ell \subseteq L$ is another inclusion of fields, and $i \colon K \rightarrow L$ and $j \colon k \rightarrow \ell$ are field homomorphisms such that $i|_k$ is equal to the composition of $j$ and the inclusion of $\ell$ into $L$, then, for any $n \geq 0$, the restriction homomorphism $r_{i|j} \colon \Omega_{K|k}^n \rightarrow \Omega_{L|\ell}^n$ maps $\nu(n)_{K|k}$ to $\nu(n)_{L|\ell}$ (this follows readily from the fact that the restriction homomorphisms commute with the differentials of the de Rham complex). In the special case where $k$ is the prime subfield of $K$, we simply write $\wp_K$ and $\nu(n)_K$ instead of $\wp_{K|k}$ and $\nu(n)_{K|k}$, respectively.

\subsection{The differential symbol} \label{Differentialsymbol} Let $k \subseteq K$ be an inclusion of fields. Since $da \wedge d(1-a) = -da \wedge da = 0$ in $\Omega_{K|k}^2$ for any $a \in K^*$, the assignment $\lpuresymbol{a_1}{a_n} \mapsto \llogsymbol{a_1}{a_n}$ gives rise to a well-defined $\mathbb{Z}$-graded ring homomorphism $K_*^M(K) \rightarrow \Omega_{K|k}^*$. In the case where $\mathrm{char}(k) = p>0$, this map descends, for any $n\geq 0$, to a homomorphism of abelian groups
\begin{equation*} h_{K|k} \colon k_{n,p}(K) \rightarrow \Omega_{K|k}^n \end{equation*}
known as the \emph{differential symbol}. In the case where $k$ is the prime subfield of $K$, we simply write $h_K$ instead of $h_{K|k}$. It is clear from their definition that the differential symbols commute with the restriction homomorphisms induced by compatible pairs of field homomorphisms. A slightly more subtle result is the following assertion concerning their compatibility with finite separable corestriction homomorphisms:

\begin{lemma} \label{LEMCompatibilityoftransfer} Let $k \subseteq K \subseteq L$ be a tower of fields, and let $\ell$ be a subfield of $L$ containing $k$. Suppose that $L$ is separable and finite over $K$ and that $\ell$ is separable and algebraic over $k$. Then, for any $n \geq 0$, the diagram
\begin{equation*} \xymatrixcolsep{3pc} \xymatrixrowsep{3pc} \xymatrix{
             k_{n,p}(L) \ar[d]_{c_{K \subseteq L}} \ar[r]^{h_{L|\ell}} & \Omega_{L|\ell}^n \ar[d]^{c_{K \subseteq L|k \subseteq \ell}} \\
             k_{n,p}(K) \ar[r]^{h_{K|k}} & \Omega_{K|k}^n} \end{equation*}
is commutative. 
\begin{proof} Let $k_0$ denote the prime subfield of $k$. By the definition of the differential symbol, the extended diagram
\begin{equation} \label{eq2.1} \xymatrixcolsep{3pc} \xymatrixrowsep{3pc} \xymatrix{
             k_{n,p}(L) \ar[d]_{c_{K \subseteq L}} \ar[r]^{h_L} & \Omega_L^n \ar[r]^{r_{|k_0 \subseteq \ell}} \ar[d]_{c_{K \subseteq L}} & \Omega_{L|\ell}^n \ar[d]^{c_{K \subseteq L|k \subseteq \ell}} \\
             k_{n,p}(K) \ar[r]^{h_K} & \Omega_K^n \ar[r]^{r_{|k_0 \subseteq k}} & \Omega_{K|k}^n} \end{equation}
includes the diagram of interest as its outer square. By the definition of the finite separable corestriction homomorphisms, the right square of \eqref{eq2.1} is commutative. In order to prove the lemma, we may therefore assume that $\ell = k = k_0$. A proof of the desired assertion this case may be found in \cite[Lemma 9.5.4]{GilleSzamuelyGaloisCohomology}. \end{proof} \end{lemma}

Now, for any fixed base field $F$ of characteristic $p>0$, and any $n \geq 0$, the maps $h_{K}$ collectively give rise to a morphism of functors
\begin{equation} \label{differentialsymbol} h \colon k_{n,p} \rightarrow \Omega^n, \end{equation}
from the category $\textsf{\textsl{Fields}}_F$ (of fields containing $F$) to the category of abelian groups. The following important result is due independently to Bloch-Kato (\cite{BlochKatopadicetalecohomology}) and Gabber (unpublished):

\begin{theorem}[{Bloch-Gabber-Kato, cf. \cite[Theorem 2.1]{BlochKatopadicetalecohomology}}] \label{BlochGabberKatoTheorem} For any $n \geq 0$, the map \eqref{differentialsymbol} is an injective morphism of functors. \end{theorem}

As an immediate corollary, we obtain:

\begin{corollary} \label{CORidentification} Let $F \subseteq L$ be an inclusion of fields of characteristic $p>0$. Then, for any $n \geq 0$, the map $h_F \colon k_{n,p}(F) \rightarrow \Omega_F^n$ induces an isomorphism of abelian groups $k_{n,p}(L/F) \simeq h_F\big(k_{n,p}(F)\big)\!\cap \Omega^n(L/F)$. \end{corollary}

Note here that, in view of \eqref{ActionofArtinSchreier}, the image $h_F\big(k_{n,p}(F)\big)$ of the differential symbol $h_F$ lies, for any $n \geq 0$, in the kernel $\nu(n)_F$ of the Artin-Schreier homomorphism $\wp_F$ of \eqref{ArtinSchreier}. In \cite{KatoGaloisCohomologyCompleteDiscrete}, Kato showed that, in fact, equality holds here: we have $h_F\big(k_{n,p}(F)\big) = \nu(n)_F$ for all $n \geq 0$. Although we shall make use of its proof, we will not explicitly appeal to this fact below. We will, however, need the following relative version of the $n=1$ case, which was proved rather earlier by Cartier (at least in the absolute case):

\begin{theorem}[{Cartier, see \cite[Ch. 2, \S 6]{CartierQuestionsRationalite}}] \label{THMCartier} Let $k \subseteq K$ be an inclusion of fields of characteristic $p>0$. Then $h_{K|k}\big(k_{n,p}(K)\big) = \nu(1)_{K|k}$. 
\begin{proof} See \cite[Theorem 9.3.3]{GilleSzamuelyGaloisCohomology}. \end{proof} \end{theorem}

\section{The norm field of a nowhere-smooth hypersurface} \label{NormField}

Let $X$ be a variety over a field $F$. We say that $X$ is \emph{generically smooth} if the locus of smooth points on $X$ is non-empty (i.e. if, for every overfield $L$ of $F$, the scheme $X \otimes_F L$ admits a regular point). Since the smooth locus of any scheme is Zariski-open, this simply amounts to asking that $X$ be smooth at its generic point. By standard results on regularity and smoothness, the latter condition may be reformulated in terms of the space of differential 1-forms of the function field $F(X)$: $X$ is smooth at its generic point if and only if $\mathrm{dim}_{F(X)}\Omega_{F(X)|F}^1 = \mathrm{dim}\;X$ (see \cite[(17.15.5)]{EGAIVPart4}). In view of Propositions \ref{PROPseparability} and \ref{PROPseparabilityfinitelygenerated}, we therefore have the following alternative characterisation of generic smoothness:

\begin{proposition} \label{PROPgenericsmoothness} Let $X$ be a variety over a field $F$. Then the following are equivalent:
\begin{enumerate} \item $X$ is generically smooth.
\item $F(X)$ is separable over $F$.
\item For any subfield $k \subseteq F$ and any $n \geq 0$, the $F$-linear homomorphism $\alpha \colon \Omega_{F|k}^n \otimes_F F(X) \rightarrow \Omega_{F(X)|k}^n$ induced by $r_{F \subseteq F(X)}$ is injective. \end{enumerate} \end{proposition}

\begin{remark} The conditions of Proposition \ref{PROPgenericsmoothness} are also equivalent to the \emph{geometric reducedness} of $X$ (i.e. the reducedness of $X \otimes_F \overline{F}$, see \cite[(4.6.1)]{EGAIVPart2}). As a simple matter of personal preference, we choose to work below with the notion of generic smoothness. \end{remark}

\begin{examples} \label{Examplesgenericallysmooth}Let $F$ be a field.
\begin{enumerate}[leftmargin=*] \item If $F$ is perfect, then every $F$-variety $X$ is generically smooth. Indeed, in this case, every finitely generated extension of $F$ is separable over $F$, and so the assertion follows immediately from Proposition \ref{PROPgenericsmoothness}.
\item Suppose that $\mathrm{char}(F) = p>0$, and let $f \in F[T_1, \hdots, T_n]$ be an irreducible polynomial in $n$ algebraically-independent variables $T_1, \hdots, T_n$ over $F$. If $f \notin F[T_1^p,\hdots,T_n^p]$, then the affine hypersurface $X_f = \lbrace f = 0 \rbrace \subset \mathbb{A}_F^n$ is generically smooth. Indeed, after reordering the $T_i$ if necessary, we may assume that $f \notin F[T_1^p,T_2,\hdots,T_n]$. Letting $L = F(T_2,\hdots,T_n)$, we then see that the function field $F(X_f)$ is $F$-isomorphic to $L(\alpha)$ for some element $\alpha$ belonging to the separable closure of $L$. The field $F(X_f)$ is therefore separably generated (and hence separable, see \cite[Theorem 26.2]{MatsumuraCommutativeRingTheory}) over $F$, and we again conclude using Proposition \ref{PROPgenericsmoothness}. \end{enumerate} \end{examples}

In the case where $X$ fails to be generically smooth, we will say that $X$ is \emph{nowhere smooth}. By definition, this means that there exists an overfield $L$ of $F$ such that the regular locus of the scheme $X \otimes_F L$ is empty. The purpose of this section is to highlight the following refinement of the latter condition in the special case where $X$ is (birationally isomorphic to) an affine or projective $F$-hypersurface:

\begin{proposition} \label{PROPExistenceofnormfield} Let $X$ be a nowhere-smooth integral affine or projective hypersurface over a field $F$. Then there exists a smallest overfield $N(X)$ of $F$ such that the regular locus of the scheme $X \otimes_F N(X)$ is empty. \end{proposition}

The field $N(X)$ of Proposition \ref{PROPExistenceofnormfield} will be called the \emph{norm field of $X$} (this terminology comes from \cite{HoffmannDiagonalFormsofDegreepinCharp},\cite{DolphinHoffmannDifferentialandBilinearForms} where $N(X)$ was defined in a more explicit way, see Remark \ref{Normfieldremark} below). Since the regular locus of any scheme is Zariski-open and dense (see \cite[(6.12.5)]{EGAIVPart2}), $N(X)$ is an invariant of the birational class of $X$. As the proof of Proposition \ref{PROPExistenceofnormfield} (given below) shows, it is a finite height-one purely inseparable extension of $F$ which is readily determined in terms of the coefficients of the polynomial (resp. homogenous polynomial) defining $X$ as a hypersurface in the ambient affine (resp. projective) space. In order to prove Proposition \ref{PROPExistenceofnormfield}, we will appeal to the following well-known statement:

\begin{lemma} \label{LEMnowhereregular} Let $X$ be a scheme. Then the regular locus of $X$ is empty if and only if, for every generic point $x$ of $X$, the Artin local ring $\mathcal{O}_{X,x}$ is non-reduced. 
\begin{proof} Since the regular locus of $X$ is Zariski-open and dense (\cite[(6.12.5)]{EGAIVPart2}), this simply amounts to the (obvious) fact that an Artin local ring is regular if and only if it is reduced (i.e. a field).\end{proof} \end{lemma}

Let us now make this more explicit in the special case where $X$ is an affine hypersurface.

\begin{lemma} \label{LEMnowhereregularhypersurface} Let $f \in F[T_1,\hdots,T_n]$ be a non-constant polynomial in $n$ algebraically-independent variables $T_1,\hdots,T_n$ over a field $F$, and let $X_f = \lbrace f = 0 \rbrace \subset \mathbb{A}^n$ be the affine $F$-hypersurface defined by its vanishing. Then the regular locus of $X_f$ is empty if and only if every irreducible divisor of $f$ in $F[T_1,\hdots,T_n]$ has multiplicity $>1$. \end{lemma}

\begin{remark} By the \emph{multiplicity} of an irreducible divisor $g$ of $f$, we mean the largest positive integer $r$ such that $g^r$ divides $f$ in $F[T_1,\hdots,T_n]$. \end{remark}

\begin{proof} If $g$ is an irreducible divisor of $f$ in $F[T_1,\hdots,T_n]$, then the affine $F$-hypersurface $X_g = \lbrace g = 0 \rbrace \subset \mathbb{A}_F^n$ is an irreducible component of $X_f$. If $x$ is the generic point of $X_g$, then the local ring $\mathcal{O}_{X_f,x}$ is $K$-isomorphic to the ring $R = F[T_1,\hdots,T_n,g^{-1}]/(g^r)$, where $r$ denotes the multiplicity of $g$ in $f$. If $r>1$, then (the class of) $g$ is a non-zero nilpotent element of $R$. On the other hand, if $r =1$, then $R$ is the fraction field of the integral domain $F[T_1,\hdots,T_n]/(g)$. In view of Lemma \ref{LEMnowhereregular}, the result follows. \end{proof}

In light of Example \ref{Examplesgenericallysmooth} (1), nowhere-smooth varieties exist only over imperfect fields (necessarily of positive characteristic). We now specialise to this situation:

\begin{proposition} \label{PROPnowheresmoothirreducible} Let $F$ be an imperfect field of characteristic $p>0$, and let $X_f \subset \mathbb{A}_F^n$ be the integral affine $F$-hypersurface defined by the vanishing of an irreducible polynomial $f \in F[T_1,\hdots,T_n]$ in $n$ algebraically-independent variables $T_1,\hdots,T_n$. Then:
\begin{enumerate} \item For any overfield $L$ of $F$, the regular locus of the scheme $X_f \otimes_F L$ is empty if and only if $f \otimes_F L = ah^p$ for some $h \in L[T_1,\hdots,T_n]$ and $a \in L^*$.
\item $X_f$ is nowhere smooth if and only if $f \in F[T_1^p,\hdots,T_n^p]$. \end{enumerate}
\begin{proof} (1) Let $L$ be an overfield of $F$. In view of Lemma \ref{LEMnowhereregularhypersurface}, we must show that $f \otimes_F L$ is a scalar multiple of a $p^{th}$ power in the ring $L[T_1,\hdots,T_n]$ if and only if all of its irreducible divisors in $L[T_1,\hdots,T_n]$ have multiplicity $>1$. The left-to-right implication here is trivial. Conversely, suppose that every irreducible divisor of $f \otimes_F L$ in $L[T_1,\hdots,T_n]$ has multiplicity $>1$. Let $g$ be any such divisor, and let $r$ denote its multiplicity. In order to show that $f \otimes_F L = ah^p$ for some $h \in L[T_1,\hdots,T_n]$ and $a \in L^*$, it will be enough to show that $r$ is a power of $p$. Re-ordering the $T_i$ if neccessary, we may assume that $f$ has non-zero degree in the variable $T_1$. Let $K = F(T_2,\hdots,T_n)$, let $M = L(T_2,\hdots,T_n)$, and let $\overline{f}$ (resp. $\overline{g}$) denote the image of $f$ (resp. $g$) under the canonical inclusion $F[T_1,\hdots,T_n] \subset K[T_1]$ (resp. $L[T_1,\hdots,T_n] \subset M[T_1]$). By Gauss' Lemma, $\overline{f}$ and $\overline{g}$ are irreducible in $K[T_1]$ and $M[T_1]$, respectively. Furthermore, as an irreducible divisor of $\overline{f} \otimes_K M$,  $\overline{g}$ has multiplicity $r$. But since $\overline{f}$ is irreducible, it then follows from elementary field theory that $r$ is a power of $p$ (e.g. \cite[Ch. V, \S 6, Prop. 1]{LangAlgebra}). This proves (1).

(2) The left-to-right implication follows from Example \ref{Examplesgenericallysmooth} (2). Conversely, if $f \in F[T_1^p,\hdots,T^p]$, then $f \otimes_F \overline{F}$ is a $p^{th}$ power in the ring $\overline{F}[T_1,\hdots,T_n]$, and so the regular locus of $X_f \otimes_F \overline{F}$ is empty by (1). In other words, $X_f$ is nowhere smooth, as we wanted.\end{proof} \end{proposition}

We can now prove Proposition \ref{PROPExistenceofnormfield}.

\begin{proof}[Proof of Proposition \ref{PROPExistenceofnormfield}] By Example \ref{Examplesgenericallysmooth} (1), $F$ is imperfect of prime characteristic $p$, say. Since the regular locus of any scheme is Zariski-open and dense (\cite[(6.12.5)]{EGAIVPart2}), we are free to replace $X$ with any variety belonging to the same birational class. In particular, we may assume that $X = X_f = \lbrace f = 0 \rbrace \subset \mathbb{A}_F^n$ for some irreducible polynomial $f \in F[T_1,\hdots,T_n]$ in $n$ algebraically-independent variables $T_1,\hdots,T_n$. By Proposition \ref{PROPnowheresmoothirreducible} (2), we necessarily have $f \in F[T_1^p,\hdots,T_n^p]$. Let $C(f) \subset F^*$ denote the set of all (non-zero) coefficients of $f$, and let $N(X) = F\big(\sqrt[p]{a/b}\;|\;a,b \in C(f)\big)$. Then $f \otimes_F N(X)$ is a scalar multiple of a $p^{th}$ power in the ring $N(X)[T_1,\hdots,T_n]$, and so the regular locus of the scheme $X \otimes_F N(X)$ is empty by Proposition \ref{PROPnowheresmoothirreducible} (1). On the other hand, if $L$ is any overfield of $F$ such that the regular locus of the scheme $X \otimes_F L$ is empty, then the same result shows that $L$ necessarily contains $N(X)$. In other words, $N(X)$ is the smallest overfield of $F$ such that the regular locus of the scheme $X \otimes_F N(X)$ is empty. \end{proof}

\begin{remark} \label{Normfieldremark}The norm field of a nowhere-smooth hypersurface $X$ over an imperfect field $F$ of characteristic $p>0$ was formally introduced by Dolphin and Hoffmann in \cite{DolphinHoffmannDifferentialandBilinearForms} (see also \cite{HoffmannDiagonalFormsofDegreepinCharp} for an earlier appearance of the same invariant). The invariant defined in \emph{loc. cit.} is actually the subfield of $p^{th}$ powers in our invariant $N(X)$ (though any distinction to be drawn between the two is clearly immaterial). Moreover,  Dolphin and Hoffmann's definition corresponds to the explicit one given in the proof of Proposition \ref{PROPExistenceofnormfield} above (in terms of the coefficients of the polynomial defining $X$ as a hypersurface in the ambient space; in particular, our interpretation of $N(X)$ as the smallest overfield of $F$ such that the regular locus of $X \otimes_F N(X)$ is empty does not appear in \emph{loc. cit.}). Nevertheless, despite this seemingly less intrinsic definition, the birational invariance of $N(X)$ may also be seen from \cite[Theorem 8.5]{DolphinHoffmannDifferentialandBilinearForms} (see also Theorem \ref{DolphinHoffmanntheorem} below).\end{remark}

\section{K\"{a}hler differentials and function fields of hypersurfaces} \label{Kahlerdifferentialsandfunctionfieldsofhypersurfaces}

Let $X$ be a variety over a field $F$. In this section we recall the explicit computation of the groups $\Omega^n(F(X)/F)$ given by Dolphin and Hoffmann in \cite{DolphinHoffmannDifferentialandBilinearForms}. We begin by noting the following immediate consequence of Proposition \ref{PROPgenericsmoothness}:

\begin{lemma} \label{LEMSeparabletriviality} Let $X$ be a variety over a field $F$. If $X$ is generically smooth, then $\Omega^n(F(X)/F) = 0$ for all $n \geq 0$. \end{lemma}

We are thus left to consider the nowhere-smooth case, which, in view of Example \ref{Examplesgenericallysmooth} (1), can only arise when $F$ is imperfect. The main result of \cite{DolphinHoffmannDifferentialandBilinearForms} is the following theorem:

\begin{theorem}[{Dolphin-Hoffmann, \cite[Theorem 8.5]{DolphinHoffmannDifferentialandBilinearForms}}] \label{DolphinHoffmanntheorem} Let $X$ be a nowhere-smooth integral affine or projective hypersurface over an imperfect field $F$ of characteristic $p>0$. Then, for any $n \geq 0$, we have $\Omega^n(F(X)/F) = \lbrace \omega \in \Omega_F^n\;|\;\omega \wedge da = 0 \text{ for all $a \in N(X)$} \rbrace$. \end{theorem}

\begin{remark} \label{Explicitremark} Let us further clarify the statement of Theorem \ref{DolphinHoffmanntheorem}. Suppose that $X$ is defined as a hypersurface in the ambient affine (resp. projective) space by the vanishing of an irreducible polynomial (resp. homogeneous polynomial) $f \in F[T_1,\hdots,T_n]$ in $n$ algebraically-independent variables $T_1,\hdots,T_n$ over $F$. Let $C(f) \subset F^*$ denote the set of all (non-zero) coefficients of $f$. Then, by the proof of Proposition \ref{PROPExistenceofnormfield}, we have $N(X) = F\big(\sqrt[p]{a/b}\;|\;a,b\in C(f) \big)$ (note, in particular, that $N(X)$ is finite over $F$). Let $b_1,\hdots,b_m$ be any $p$-basis of $N(X)^p$ over $F^p$ (again, this is a finite set). Since $N(X)$ is a height-one purely inseparable extension of $F$, the $b_i$ all belong to $F^*$. Theorem \ref{DolphinHoffmanntheorem} then asserts that, for any $n \geq 0$, we have $\Omega^n(F(X)/F) = \Omega_F^{n-m} \wedge \big(\llogsymbol{b_1}{b_m}\big)$. In particular, we see that $\Omega^n(F(X)/F) = 0$ whenever $n < \mathrm{log}_p[N(X):F]$. \end{remark}

In view of Remark \ref{Explicitremark}, Theorem \ref{DolphinHoffmanntheorem} therefore yields:

\begin{corollary} \label{CORconjectureimplication} A positive solution to Conjecture \ref{Secondaryconjecture} (in degree $n$) implies a positive solution to Conjecture \ref{Mainconjecture} (in degree $n$). \end{corollary}

\section{Proof of main results}

We now prove the main results of this article, beginning with:

\begin{proposition} \label{PROPproofofsecondaryconjecture} Conjecture \ref{Secondaryconjecture} is true for $n=m$. \end{proposition}

The proof of Proposition \ref{PROPproofofsecondaryconjecture} is, in fact, essentially implicit in that of \cite[Theorem 1]{KatoGaloisCohomologyCompleteDiscrete}. For lack of a precise reference, however, we provide full details here. As in \emph{loc. cit.} we will need the following lemma:

\begin{lemma}[{cf. \cite[Lemma 3]{KatoGaloisCohomologyCompleteDiscrete}}] \label{LEMmaintechnicallemma} Let $k \subseteq K$ be an inclusion of fields with $K$ purely inseparable and finite of degree $p$ over $k$, and let $g \colon K \rightarrow k$ be a $k$-linear homomorphism. Then there exist an overfield $\ell$ of $k$ and a non-zero element $c \in L \coloneqq \ell \cdot K$ such that:
\begin{enumerate} \item $\ell$ is finite of prime-to-$p$ degree over $k$.
\item $\overline{g}(c^i) = 0$ for all $1 \leq i <p$, where $\overline{g}$ denotes the unique $\ell$-linear homomorphism $L \rightarrow \ell$ sastifying $\overline{g}(a) = g(a)$ for all $ a \in K$. \end{enumerate} \end{lemma}

The proof of Lemma \ref{LEMmaintechnicallemma} begins with the following auxiliary result:

\begin{lemma} \label{LEMsecondtechnicallemma} Let $k \subseteq K$ be an inclusion of fields with $K$ with inseparable and finite of degree $p$ over $k$, and let $\omega \in \Omega_{K|k}^1 \setminus B_{K|k}^1$. Then there exist an overfield $\ell$ of $k$ and non-zero elements $u \in \ell$ and $y \in L \coloneqq \ell \cdot K$ such that:
\begin{enumerate} \item $\ell$ is finite of prime-to-$p$ degree over $k$.
\item $r_{K \subseteq L|k \subseteq \ell}(\omega) = u\logsymbol{y}$ in $\Omega_{L|\ell}^1$. \end{enumerate}
\begin{proof} Let $b \in K$ be such that $K = k(b)$. The $K$-vector space $\Omega_{K|k}^1$ is one dimensional, generated by the element $\logsymbol{b}$. In particular, there exists $a \in K^*$ such that $\omega = a\logsymbol{b}$ in $\Omega_{K|k}^1$. On the other hand, Proposition \ref{PROPCohomologydeRhamcharp} shows that the space $H^1(\Omega_{K|k}^\bullet) = \Omega_{K|k}^1/B_{K|k}^1$ is one dimensional over $k$. Since $\omega \notin B_{K|k}^1$, it follows that there exists $\rho \in k^*$ such that
\begin{equation} \label{eq5.1} a^p\logsymbol{b} \equiv \rho \omega \equiv \rho a\logsymbol{b} \pmod{B_{K|k}^1}. \end{equation}
Let $u$ be any $(p-1)^{\mathrm{st}}$ root of $\rho$ in $\overline{k}$, and let $\ell = k(u)$ and $L = \ell \cdot K = K(u)$. Clearly $\ell$ has prime-to-$p$ degree over $k$. Since restriction homomorphisms on differential forms are compatible with the differentials of the de Rham complex, applying the map $r_{K \subseteq L| k \subseteq \ell}$ to \eqref{eq5.1} yields the congruence
\begin{equation} \label{eq5.2} a^p \logsymbol{b} \equiv \rho a \logsymbol{b} \pmod{B_{L|\ell}^1} \end{equation}
in $\Omega_{L|\ell}^1$. Dividing both sides of \eqref{eq5.2} by $u^p \in \ell^*$ and rearranging, we see that
\begin{equation*} (c^p - c) \logsymbol{b} \in B_{L|\ell}^1, \end{equation*}
where $c = a/u \in L^*$. In other words (see \eqref{ActionofArtinSchreier}), we have $c\logsymbol{b} \in \nu(1)_{L|\ell}$. By Theorem \ref{THMCartier}, it follows that there exists $y \in L^*$ such that $c \logsymbol{b} = \logsymbol{y}$ in $\Omega_{L|\ell}^1$. Since $r_{K \subseteq L|k \subseteq \ell}(\omega) = uc\logsymbol{b}$ in $\Omega_{L|\ell}^1$, this proves the desired assertion. \end{proof} \end{lemma}

We also need:

\begin{lemma} \label{LEMthirdbasiclemma} Let $k \subseteq K$ be an inclusion of fields with $K$ finite over $k$, and let $V$ and $W$ be codimension-one $k$-linear subspaces of $K$. Then there exists $\alpha \in K^*$ such that $\alpha \cdot V = W$.
\begin{proof} We may assume that $K \neq k$. Let $n = [K:k]>1$, and let $v_1,\hdots,v_{n-1}$ be a basis of $V$ over $k$. Let $\alpha$ be any non-zero element of $K$. Then $\alpha \cdot V = W$ if and only if $\alpha v_i \in W$ for all $1 \leq i <n$. For each such $i$, let $W_i = (v_i)^{-1} \cdot W \subset K$. Then the $W_i$ are codimension-one $k$-linear subspaces of $K$, and we have $\alpha \cdot V = W$ if and only if $\alpha$ lies in the intersection of the $W_i$. In order to prove the lemma, it therefore suffices to check that the $W_i$ intersect non-trivially. This is clear for dimension reasons, and so the statement follows. \end{proof} \end{lemma}

Now we can prove Lemma \ref{LEMmaintechnicallemma}:

\begin{proof}[Proof of Lemma \ref{LEMmaintechnicallemma}] Let $b \in K$ be such that $K = k(b)$. The $K$-vector space $\Omega_{K|k}^1$ is one dimensional, generated by the element $\logsymbol{b}$. Let $\phi \colon \Omega_{K|k}^1 \rightarrow K$ be the unique $K$-linear isomorphism which sends the latter element to $1$, and let $g' = g \circ \phi \colon \Omega_{K|k}^1 \rightarrow k$. Now, if $g(1) = 0$, then the lemma holds trivially with $\ell = k$ and $c =1$. We may therefore assume that $g(1) \neq 0$. In this case, the map $g'$ is non-trivial, and so $\mathrm{ker}(g')$ is a codimension-one $k$-linear subspace of $\Omega_{K|k}^1$. At the same time, Proposition \ref{PROPCohomologydeRhamcharp} shows that
\begin{equation} \label{eq5.3} B_{K|k}^1 = \bigoplus_{i=1}^{p-1} k \cdot b^i \logsymbol{b} \end{equation}
is also a codimension-one $k$-linear subspace of $\Omega_{K|k}^1$. In view of Lemma \ref{LEMthirdbasiclemma}, it follows that there exists an element $\alpha \in K^*$ such that $B_{K|k}^1 = \alpha \cdot \mathrm{ker}(g')$ in $\Omega_{K|k}^1$. If $a$ is any element of $K$, we then see that $g(a) = 0$ if and only if $\alpha a\logsymbol{b} \in B_{K|k}^1$. In particular, since $g(1) \neq 0$, the element $\omega \coloneqq \alpha \logsymbol{b} \in \Omega_{K|k}^1$ does not lie in $B_{K|k}^1$. Lemma \ref{LEMsecondtechnicallemma} therefore implies the existence of an overfield $\ell$ of $k$ and non-zero elements $u \in \ell$ and $y \in L \coloneqq \ell \cdot K$ such that
\begin{enumerate} \item $\ell$ is finite of prime-to-$p$ degree over $k$, and
\item $r_{K \subseteq L|k \subseteq \ell}(\omega) = u\logsymbol{y}$ in $\Omega_{L|\ell}^1$. \end{enumerate}
Let $\overline{g} \colon L \rightarrow \ell$ be the unique $\ell$-linear homomorphism satisfying $\overline{g}(a) = g(a)$ for all $a \in K$. If we can show that $\overline{g}(y^i) = 0$ for all $1 \leq i <p$, then we can take $c = y$ and the lemma will be proved. First, note that since $\ell$ is separable and finite over $k$, the natural $L$-linear homomorphism $\beta \colon \Omega_{K|k}^1 \otimes_K L \rightarrow \Omega_{L|\ell}^1$ induced by $r_{K \subseteq L|k \subseteq \ell}$ is an isomorphism by Proposition \ref{PROPseparablealgebraic}. Moreover, in view of \eqref{eq5.3}, this map $\beta$ restricts to an isomorphism $B_{K|k}^1 \otimes_K L \xrightarrow{\sim} B_{L|\ell}^1$. Now, let $\overline{g}' \colon \Omega_{L|\ell}^1 \rightarrow L$ denote the $\ell$-linear composition $(g' \otimes \overline{g}) \circ \beta^{-1}$. Then, by construction, we have $\overline{g}(a) = \overline{g}'(a\logsymbol{b})$ for all $a \in L$. Moreover, $\beta$ induces an isomorphism $\mathrm{ker}(g') \otimes_k \ell \simeq \mathrm{ker}(\overline{g}')$. It therefore follows that $B_{L|\ell}^1 = \alpha \cdot \mathrm{ker}(\overline{g}')$ in $\Omega_{L|\ell}^1$. In particular, for any $a \in L$, we see that $\overline{g}(a) = 0$ if and only if $\alpha a \logsymbol{b} \in B_{L|\ell}^1$. Finally, let $1 \leq i < p$, and let $1 \leq j <p$ be such that $ij = 1\pmod{p}$. Then
\begin{equation*} \alpha y^i \logsymbol{b} = uy^{i-1}dy = jud(y^i) = d(juy^i) \in B_{L|\ell}^1, \end{equation*}
and so $\overline{g}(y^i) = 0$. Since this holds for all such $i$, the lemma is proved. \end{proof}

Next, we need the following reduction statement:

\begin{lemma} \label{LEMtransferlemma} Let $k \subseteq K \subseteq L$ be a tower of fields, and let $\ell$ be a subfield of $L$ containing $k$. Suppose that $L$ is finite of prime-to-$p$ degree over $K$, and that $\ell$ is separable and algebraic over $k$. Let $\omega \in \Omega_{K|k}^n$ for some $n \geq 0$. Then $\omega \in h_{K|k}\big(k_{n,p}(K)\big)$ if and only if $r_{K \subseteq L|k \subseteq \ell}(\omega) \in h_{L|\ell}\big(k_{n,p}(L)\big)$.
\begin{proof} Since the differential symbols are compatible with restriction homomorphisms, the left-to-right implication is clear. Conversely, if $r_{K \subseteq L|k \subseteq \ell}(\omega) \in h_{L|\ell}\big(k_{n,p}(L)\big)$, then it follows from Lemma \ref{LEMCompatibilityoftransfer} and the projection formula (see \S \ref{Kahlerdifferentials} above) that $[L:K]\omega \in h_{L|\ell}\big(k_{n,p}(K)\big)$. Since $[L:K]$ is invertible modulo $p$ by assumption, the result follows. \end{proof} \end{lemma}

Finally, we need the following basic observation:

\begin{lemma} \label{LEMreductionlemma} Let $k \subseteq K$ be an inclusion of fields of characteristic $p>0$ such that $K^p \subseteq k$. Let $b_1,\hdots,b_n \in K$ be such that $[k(b_1,\hdots,b_n):k] = p^n$, and let $a \in K$ be such that $a \llogsymbol{b_1}{b_n} \in \nu(n)_{K|k}$. Then $a \in k(b_1,\hdots,b_n)$ and $a\llogsymbol{b_1}{b_n} \in \nu(n)_{k(b_1,\hdots,b_n)|k}$.
\begin{proof} Let $k' = k(b_1,\hdots,b_n)$. Since $a\llogsymbol{b_1}{b_n} \in \nu(n)_{K|k}$, we have
\begin{equation} \label{eq5.4} (a^p-a)\llogsymbol{b_1}{b_n} \in B_{K|k}^n \end{equation}
(see \eqref{ActionofArtinSchreier} above). Applying the differential $d \colon \Omega_{K|k}^n \rightarrow \Omega_{K|k}^{n+1}$ to both sides of \eqref{eq5.4}, we see that $da \wedge \llogsymbol{b_1}{b_n} = 0$ in $\Omega_{K|k}^{n+1}$, from which it immediately follows that $a \in k'$. To prove the second statement, we may assume that $K$ is finite over $k$, say of degree $p^{n+r}$ for some $r \geq 0$. If $r = 0$, then $K = k'$, and there is nothing to prove. Otherwise, we can extend $b_1,\hdots,b_n$ to a $p$-basis $b_1,\hdots,b_n,b_{n+1},\hdots,b_{n+r}$ of $K$ over $k$. Since $a^p - a \in k'$, a quick inspection of the additive decomposition of the complex $\Omega_{K|k}^\bullet$ associated to this choice of $p$-basis (as per \eqref{additivedecompositiondeRham} above) shows that \eqref{eq5.4} already holds over $k'$. In other words, we have
\begin{equation*} (a^p - a) \llogsymbol{b_1}{b_n} \in B_{K|k}^1. \end{equation*}
This means that $a \llogsymbol{b_1}{b_n} \in \nu(n)_{k'|k}$ (again, see \eqref{ActionofArtinSchreier}), as we wanted. \end{proof} \end{lemma}

Now, Proposition \ref{PROPproofofsecondaryconjecture} will follow readily from the following result:

\begin{proposition} \label{PROPkeyresult} Let $k \subseteq K$ be an inclusion of fields of characteristic $p>0$ such that $K^p \subset k$. Let $b_1,\hdots,b_n \in K$ be such that $[k(b_1,\hdots,b_n):k] = p^n$, and let $a \in K$ be such that $\omega \coloneqq a \llogsymbol{b_1}{b_n} \in \nu(n)_{K|k}$. Then $\omega$ may be expressed in $\Omega_{K|k}^n$ as a sum of elements of the form $\llogsymbol{a_1}{a_n}$, where $k(a_1,\hdots,a_n) = k(b_1,\hdots,b_n)$.
\begin{proof} By Lemma \ref{LEMreductionlemma}, we may assume that $K = k(b_1,\hdots,b_n)$. In this case, proving the lemma amounts to showing that $\omega \in h_{K|k}\big(k_{n,p}(K)\big)$, and we proceed by induction on $n$. Now, the $K$-vector space $\Omega_{K|k}^n$ is one dimensional, generated by the element $\eta \coloneqq \llogsymbol{b_1}{b_n}$. By assumption, we have
\begin{equation} \label{eq5.5} a^p\eta - \omega = (a^p-a)\eta \in B_{K|k}^n \end{equation}
(see \eqref{ActionofArtinSchreier} above). Let $k_1 = k(b_1)$, and consider the $k$-linear map $g \colon k_1 \rightarrow H^n(\Omega_{K|k}^\bullet) = \Omega_{K|k}^n/B_{K|k}^n$ defined by the assignment $\alpha \mapsto \alpha \omega \pmod{B_{K|k}^n}$. By Proposition \ref{PROPCohomologydeRhamcharp}, $H^n(\Omega_{K|k}^\bullet)$ is one dimensional over $k$, generated by $[\eta]$, and may therefore be identified with $k$ as a $k$-vector space. By Lemma \ref{LEMmaintechnicallemma}, there exists an overfield $\ell$ of $k$ and a non-zero element $c \in \ell_1 \coloneqq \ell \cdot k_1$ such that
\begin{enumerate} \item $\ell$ is finite of prime-to-$p$ degree over $k$, and 
\item $\overline{g}(c^i) = 0$ for all $1 \leq i < p$, where $\overline{g}$ denotes the unique $\ell$-linear map $\ell_1 \rightarrow \ell \simeq H^n(\Omega_{\ell \cdot K|\ell}^\bullet)$ satisfying $\overline{g}(\alpha) = g(\alpha)$ for all $\alpha \in k_1$. \end{enumerate}
On the other hand, Lemma \ref{LEMtransferlemma} shows that, in order to prove the proposition, we are free to replace $k$ by $\ell$, $k_1$ by $\ell_1$, $K$ by $L \coloneqq \ell \cdot K$ and $\omega$ by $r_{K \subseteq L|k \subseteq \ell}(\omega)$ (note that the latter element belongs to $\nu(n)_{L|\ell}$ by the compatibility of restriction homomorphisms with the differentials of the de Rham complex). We may therefore assume that there exists $c \in k_1^*$ such that $g(c^i) = 0$ for all $1 \leq i <p$. In other words, we have $c \in k_1^*$ such that
\begin{equation} \label{eq5.6} c^i \omega \in B_{K|k}^n \end{equation}
for all such $i$. Now, if $c \in k$, then, since $B_{K|k}^n$ is closed under scalar multiplication by $k$, \eqref{eq5.6} (with $i=1$, say) implies that $\omega \in B_{K|k}^n$. By \eqref{eq5.5}, it then follows that $a^p\eta \in B_{K|k}^n$. But since $K^p \subseteq k$, this implies that $\eta \in B_{K|k}^n$, which is impossible in view of the fact that $[\eta]$ generates the one-dimensional $k$-vector space $H^n(\Omega_{K|k}^\bullet)$. It follows that $k_1 = k(c)$. In particular, $c,b_2,\hdots,b_n$ is a $p$-basis of $K$ over $k$. Setting
\begin{equation*} \eta' = \begin{cases} \llogsymbol{b_2}{b_n} & \text{if $n>1$} \\
                                                              1 & \text{otherwise,} \end{cases} \end{equation*}
it follows that $\Omega_{K|k}^n$ is generated over $K$ by the element $\logsymbol{c} \wedge \eta'$ (if $n=1$, this is to be understood simply as the element $\logsymbol{c}$). In particular, we can write $\omega = a'\logsymbol{c} \wedge \eta'$ for some $a' \in K^*$. In direct analogy with \eqref{eq5.5}, we have
\begin{equation} \label{eq5.7} a'^p \logsymbol{c} \wedge \eta' - \omega = (a'^p - a') \logsymbol{c} \wedge \eta' \in B_{K|k}^n. \end{equation} 
At the same time, we also have
\begin{equation} \label{eq5.8} c^ia' \logsymbol{c} \wedge \eta' \in B_{K|k}^n \end{equation}
for all $1 \leq i < p$ by \eqref{eq5.6}. Now, consider the additive decomposition
\begin{equation*} \Omega_{K|k}^\bullet = \bigoplus_{\alpha} \Omega_{K|k}^\bullet(\alpha) \end{equation*}
of the complex $\Omega_{K|k}^\bullet$ associated to the (ordered) $p$-basis $c,b_2,\hdots,b_n$ of $K$ over $k$ (as per \eqref{additivedecompositiondeRham} above). By Proposition \ref{PROPCohomologydeRhamcharp} (1), the projection of $B_{K|k}^n$ to the space $\Omega_{K|k}^n(\overline{0})$ is trivial. Applying this to the elements on the left sides of \eqref{eq5.7} and \eqref{eq5.8}, respectively, we see that
\begin{equation} \label{eq5.9} a' = a'^p - v \end{equation}
for some $v \in V$, where $V$ denotes the $k_1$-linear subspace of $K$ generated by the monomials $b_2^{\alpha_2}\cdots b_n^{\alpha_n}$ (with $0 \leq \alpha_i < p$) having at least one non-zero exponent. In particular, in the case where $n=1$ (whereby $V=0$), we have $a' \in \mathbb{Z}/p\mathbb{Z}$. Letting $m$ be any positive integer in the mod-$p$ congruence class of $a'$, we then have
\begin{equation*} \omega = a' \logsymbol{c} = h_{K|k}(\puresymbol{c^m}) \end{equation*}
in $\Omega_{K|k}^1$, and so the desired assertion holds in this case. Suppose now that $n>1$. Note that the set $b_2,\hdots,b_n$ is a $p$-basis of $K$ over $k_1$. In particular, the $K$-vector space $\Omega_{K|k_1}^{n-1}$ is one dimensional, generated by the element $r_{|k \subseteq k_1}(\eta')$. Now, in view of the definition of $V$, Proposition \ref{PROPCohomologydeRhamcharp} (2) (applied to the (ordered) $p$-basis $b_2,\hdots,b_n$ of $K$ over $k_1$) implies that $vr_{|k \subseteq k_1}(\eta') \in B_{K|k_1}^{n-1}$. By \eqref{eq5.9}, we therefore have
\begin{equation*} r_{|k \subseteq k_1}\big((a'^p-a)\eta'\big) = vr_{|k \subseteq k_1}(\eta') \in B_{K|k_1}^{n-1}. \end{equation*}
In other words, we have $r_{|k \subseteq k_1}(a'\eta') \in \nu(n-1)_{K|k_1}$ (see \eqref{ActionofArtinSchreier} above). By the induction hypothesis, it follows that
\begin{equation} \label{eq5.10} r_{|k \subseteq k_1}(a'\eta') \in h_{K|k_1}\big(k_{n-1,p}(K)\big). \end{equation}
Since any element in the kernel of $r_{|k \subset k_1}$ is divisible by $\logsymbol{c}$ \big(recall that $k_1 = k(c)$\big), \eqref{eq5.10} implies that
\begin{equation} \label{eq5.11} a'\eta' \in h_{K|k}\big(k_{n-1,p}(K)\big) + \logsymbol{c} \wedge \Omega_{K|k}^{n-2}. \end{equation}
Finally, since $\omega = a'\logsymbol{c} \wedge \eta'$, taking the wedge product of \eqref{eq5.11} with $\logsymbol{c}$ gives
\begin{equation*} \omega \in \logsymbol{c} \wedge h_{K|k}\big(k_{n-1,p}(K)\big). \end{equation*}
Since $\logsymbol{c} \wedge h_{K|k}\big(k_{n-1,p}(K)\big) \subset h_{K|k}\big(k_{n,p}(K)\big)$, this proves what we wanted. \end{proof} \end{proposition}

We can now prove Proposition \ref{PROPproofofsecondaryconjecture}:

\begin{proof}[Proof of Proposition \ref{PROPproofofsecondaryconjecture}] Since $m=n$, every element of $V$ is of the form $a \llogsymbol{b_1}{b_n}$ for some $a \in K$. Since the image of the differential symbol $h_F$ lies in $\nu(n)_F$, the desired assertion follows from Proposition \ref{PROPkeyresult} (applied to the case where $K = F$ and $k = F^p$). \end{proof}

Finally, we prove Theorem \ref{Maintheorem}:

\begin{proof}[Proof of Theorem \ref{Maintheorem}] By Corollary \ref{CORidentification}, the differential symbol $h_F$ induces an isomorphism
\begin{equation*} k_{n,p}(F(X)/F) \simeq h_F\big(k_{n,p}(F)\big)\!\cap \Omega^n(F(X)/F) \end{equation*}
for any $n \geq 0$. Parts (1) and (2) then follow immediately from Lemma \ref{LEMSeparabletriviality} and Theorem \ref{DolphinHoffmanntheorem}, respectively (in the latter case, see also Remark \ref{Explicitremark}). In the same way, part (3) follows from Theorem \ref{DolphinHoffmanntheorem} (plus Remark \ref{Explicitremark}) together with Proposition \ref{PROPproofofsecondaryconjecture}. \end{proof}

\noindent {\bf Acknowledgements.} This work was carried out as part of my Ph.D. thesis at the University of Nottingham, and I would like to thank my advisors, Detlev Hoffmann and Alexander Vishik, for very helpful discussions. In its present form, this paper was partly written during a visit to the University of California, Los Angeles. I gratefully acknowledge the support of the Cecil King Memorial Foundation and the London Mathematical Society which made this visit possible.

\providecommand{\bysame}{\leavevmode\hbox to3em{\hrulefill}\thinspace}
\providecommand{\MR}{\relax\ifhmode\unskip\space\fi MR }
\providecommand{\MRhref}[2]{%
  \href{http://www.ams.org/mathscinet-getitem?mr=#1}{#2}
}
\providecommand{\href}[2]{#2}

\end{document}